\documentclass{elsart}
\usepackage{amsmath,amscd,amssymb,nccbbb}

\theoremstyle{plain}
\newtheorem{question}[thm]{Question}

\theoremstyle{definition}

\numberwithin{equation}{section}

\newcommand\proof{\smallskip\noindent{\bf Proof.}\quad}
\newcommand\eproof{$\qed$  \\ \\}

\def \p{p}

\newcommand {\Ok}{\mathcal{O}/{\wp^k}}
\newcommand {\Gk}{\text{G}(\OO_k)}

\newcommand {\F}{\mathbf{F}}
\newcommand {\PP}{\mathbb{P}}
\newcommand {\Abb}{\mathbb{A}}

\def \al{\alpha}
\def \be{\beta}
\def \gam{\gamma}

\def \Al{a}
\def \Be{b}
\def \Gam{t}

\def \domstrict {\prec}
\def \domnonstrict {\preceq}
\def \incl {\le}

\def \arch {\textbf{\tiny{Arch}}}
\def \nonarch {\textbf{\tiny{NonArch}}}

\def \mt{\texttt{M}}

\newcommand {\Om}{\Omega}

\newcommand {\OO}{\mathcal{O}}
\newcommand {\Sc}{\mathcal{S}}

\newcommand {\FF}{\mathcal{F}}
\newcommand {\JJ}{\mathcal{J}}
\newcommand {\Ac}{\mathcal{A}}
\newcommand {\Bc}{\mathcal{B}}

\newcommand {\Uc}{\mathcal{U}}

\newcommand {\R}{\mathbb{R}}
\newcommand {\C}{\mathbb{C}}
\newcommand {\Q}{\mathbb{Q}}
\newcommand {\K}{\mathbb{K}}
\newcommand {\N}{\mathbb{N}}

\newcommand {\Kp}{{\K}}

\newcommand {\Gt}{\text{G}}
\newcommand {\Ut}{\text{U}}
\newcommand {\Ot}{\text{O}}
\newcommand {\GL}{\text{GL}}
\newcommand {\PGL}{\text{PGL}}
\newcommand {\Ker}{\text{Ker}}
\newcommand {\End}{\text{End}}
\newcommand {\Gr}{\text{Gr}}

\newcommand {\Span}{\text{Span}}

\def \HH{{\mathcal{H}}}

\def \dh {\texttt{dh}}
\def \dS {\texttt{dS}}

\def \la {\lambda}
\def \La {\Lambda}

\def \ta {\texttt{t}}

\def \f {\phi}

\def \cb {{\texttt c}}
\def \gb {{\texttt g}}
\def \eb {{\texttt e}}
\def \xb {{\mathbf x}}
\def \ub {{\mathbf u}}
\def \Ub {{\mathbf U}}

\def \Eb {{\texttt E}}
\def \Cb {{\texttt C}}

\begin{document}

\begin{frontmatter}

\title{From $p$-adic to real Grassmannians via the quantum}
\thanks{Supported by Israel Science Foundation (ISF grant no.
100146), by NWO (grant no. 613.006.573) and by Marie Curie research
training network LIEGRITS (MRTN-CT 2003-505078).}
\author{Uri Onn}
\ead{onn@math.jussieu.fr}
\address{Institut de Math\'{e}matiques de Jussieu}


\begin{abstract}
Let $\F$ be a local field. The action of $\GL_n(\F)$ on the
Grassmann variety $\Gr(m,n,\F)$ induces a continuous representation
of the maximal compact subgroup of $\GL_n(\F)$ on the space of
$L^2$-functions on $\Gr(m,n,\F)$. The irreducible constituents of
this representation are parameterized by the same underlying set
both for Archimedean and non-Archimedean fields \cite{HG2,JC}.
 This paper connects the Archimedean and non-Archimedean theories using the quantum
 Grassmannian \cite{Stokman2,Stokman1}. In particular, idempotents in the Hecke algebra
 associated to this representation are the image of the quantum zonal
 spherical functions after taking appropriate limits. Consequently, a correspondence is established
 between some irreducible representations with Archimedean and non-Archimedean origin.
\end{abstract}


\end{frontmatter}


\section{Introduction}
This paper is concerned with relationships between the Archimedean
and non-Archimedean places of a number field. Since the early works
of Weil, Artin, Iwasawa, Tate \cite{Tate} and the far reaching
conjectures of Langlands, deep relations have been discovered
between the arithmetic of a number field and the representation
theory of algebraic groups over the local fields. It is within the
framework of representation theory that the relations between the
local fields, the {\em places} of the number field, will be
discussed here.

Local fields occur naturally as the completions of global fields. A
{\em global field} is either a {\em number field}, that is a finite
extension of the rational numbers, or a {\em function field}, that
is a field of rational functions of a curve defined over a finite
field. A local field can be either {\em Archimedean} ($\R$ or $\C$)
or {\em non-Archimedean} (Laurent series over a finite field or a
finite extension of $\Q_p$).
 In the function field case all the
completions are non-Archimedean and thus carry the same nature. In
contrast, in the number field case both Archimedean and
non-Archimedean completions occur, thus having a completely
different nature. For example, the former is connected and the
latter is totally disconnected.

Let $\F$ be a local field. For a non-Archimedean field, let $\OO$ be
the ring of integers and $\wp$ be the maximal ideal. Let $K^{\F}$ be
the maximal compact subgroup of $\GL_n(\F)$, for some $n \in \N$
which will be fixed throughout this paper. We have
\begin{equation*}
K^{\F}=\left\{
         \begin{array}{ll}
           \Ot(n)=\text{the orthogonal group}, & \F=\R; \\
           \Ut(n)=\text{the unitary group},& \F=\C; \\
           \GL_n(\OO) \simeq \varprojlim \GL_n(\OO/\wp^k), & \F ~\text{non-Archimedean.}
         \end{array}
       \right.
\end{equation*}
In particular, for Archimedean fields $K^{\F}$ is a Lie group while
for non-Archimedean fields it is totally disconnected. In order to
be able to compare between them we appeal to representation theory.
In this paper we focus on a special representation of $K^{\F}$, the
{\em Grassmann representation}, which arises from its natural action
on $X_m^{\F}=\Gr(m,n,\F)$, the variety of $m$-dimensional subspaces
of a fixed $n$-dimensional space over $\F$. The natural
representation space is $L^2(X_m^{\F})$ or its dense subspace of
{\em smooth} functions $\Sc(X_m^{\F})$, with the action
\[
[g\cdot f](x)=f(g^{-1} x), \quad f \in L^2(X_m^{\F}), g \in K^{\F} .
\]
As far as the decomposition to irreducibles is concerned, there is
no difference between the two spaces. By smooth functions we mean
infinitely differentiable for Archimedean places and locally
constant for non-Archimedean ones. To define the $L^2$ structure,
the transitive action of $K^{\F}$ on $X_m^{\F}$ is used, and the
measure on $X_m^{\F}$ is taken to be the projection of the
normalized Haar measure from the group. Then, for all local fields,
for the Archimedean ones \cite{JC} and for the non-Archimedean ones
\cite{HG2}, the following decomposition holds.

\begin{thm}[James-Constantine, Hill]\label{multiplicity} For any local
field $\F$ and $m \le [\frac{n}{2}]$, the Grassmann representation
is a multiplicity free direct sum of irreducible representations of
$K^{\F}$ indexed by $\La_m$, the set of partitions with at most $m$
parts.
\end{thm}

Let $\{\mathcal{U}_{\la}^{\F}\}_{\la \in \La_m}$ be the irreducible
representations which occur in $L^2(X_m^{\F})$. In view of the
independence of the labeling set on the field, it is natural to ask
the following question.
\begin{question}\label{q1}
Fix $\la \in \La_m$. Are $\{\mathcal{U}_{\la}^{\F}\}_{\F}$ related when $\F$
runs over all local fields?
\end{question}
Our goal is to address this question. For this purpose, the {\em
Hecke algebra} of intertwining operators $\HH_m^{\F}=\Sc(X_m^{\F}
\times _{K^{\F}} X_m^{\F})$ will be used. This is the convolution
algebra of smooth functions on $X_m^{\F} \times X_m^{\F}$ which are
invariant under the diagonal action of $K^{\F}$. An element of
$\HH_m^{\F}$ defines an intertwining operator by realizing it as an
integration kernel. The measure on $\Om_m^{\F}:=X_m^{\F} \times
_{K^{\F}} X_m^{\F}$, denoted by $\dh^{\F}$, is the projection of the
Haar measure from the group and described explicitly in
\S\ref{Archimedean} and \S\ref{non-Archimedean}. As this algebra is
commutative for all local fields, the first part of Theorem
\ref{multiplicity} follows. The minimal idempotents of the algebra
have been computed in \cite{JC} for Archimedean fields and in
\cite{BO2} for non-Archimedean fields:
\begin{itemize}
\item \textbf{Archimedean fields \cite{JC}.} The minimal idempotents in the Hecke algebra are
naturally associated to polynomial representations of $\GL_m$. In
particular, they are parameterized by $\La_m$. They are
eigenfunctions of the Laplacian on the Grassmann manifold with
distinct eigenvalues.
\item \textbf{Non-Archimedean fields \cite{BO2}.} The minimal idempotents in the Hecke algebra are
naturally associated with finite quotients of $\OO^m$, the free
module of rank $m$. In particular, they are parameterized by
$\La_m$. The idempotents are computed in terms of combinatorial
invariants of the lattice of submodules of $\OO^m$.
\end{itemize}

Interestingly, geometry plays an important role in both cases;
geometrically defined operators which commute with the group action
are sufficient to separate representations. In the Archimedean case,
it is the Laplacian on the Grassmann manifold, whereas in the
non-Archimedean case, a family of discrete averaging operators plays
the same role. The identical parametrization of irreducibles is
reflected by the same parametrization of idempotents in the Hecke
algebras for the different local fields. To show the link between
the irreducibles
 labeled by the same partition for the different fields, the quantum
 Grassmannian will be used in the following scheme.

 Each of the Hecke algebras $\HH_m^{\F}$ is
 characterized by a triplet (space, measure, idempotents)
\begin{align*}
\big(&\Om_m^{\F},  ~\dh^{\F}, \{\eb_{\la}^{\F}\}_{\la \in
\La_m}\big)
\end{align*}
These will be 'interpolated' by similar objects which arise in the
quantum Grassmannian $\Ub_q(n)/\Ub_q(m) \times \Ub_q(n-m)$ (cf.
\cite{Stokman2} for a detailed discussion). The objects which will
be used are
\begin{align*}
\big(\Om^q_m, ~\dS^q_m, \{\Eb^q_{\la}\}_{\la \in \La_m}\big)
\end{align*}
The precise definition of these $q$-objects is given later on.
Roughly, $\Om^q_m$ is the $q$-exponentiation of a shift of $\La_m$;
$\dS_m^q(\xb;\Al,\Be,\Gam)$ is the  $q$-Selberg measure
\cite{Habsieger,Kadell1,Kaneko,Aomoto1} defined on $\Om^q_m$; and
$\{\Eb^q_{\la}(\xb;\Al,\Be,\Gam)\}_{\la \in \La_m}$ are the zonal
spherical functions which occur in the quantum Grassmannian.
 The zonal spherical functions, also called multivariable little $q$-Jacobi polynomials
\cite{Stokman2,Stokman1}, are orthogonal with respect to the
$q$-Selberg measure.

By taking appropriate limits, the $q$-objects interpolate between
the objects related to the local fields. In the Archimedean limit $q
\rightarrow 1$, the space $\Om^q_m$ becomes dense in the Archimedean
space, and the atomic $q$-Selberg measure approximates the
continuous Selberg measure. In the non-Archimedean limit $q
\rightarrow 0$, the space itself remains discrete, and the
$q$-measure specializes to give the non-Archimedean measure. Thus,
for any local field $\F$, the distribution $f \mapsto
\int_{\Om_m^{\F}}f\dh^{\F}$ is the limit of the distribution $f
\mapsto \int_{\Om_m^q}f\dS_m^q$ (\S\ref{interpolation}, Theorem
\ref{thm-measures}). Under the same limits the quantum zonal
spherical functions $\{\Eb^q_{\la}\}$ are mapped to
$\{\eb_{\la}^{\F}\}$ (\S\ref{interpolation}, Theorem
\ref{thm-idempotents}).

\subsection{Related works}

Similar interpolations between $p$-adic and real zonal spherical
functions using $q$-special functions have been established in
several instances. For $\PGL_2$, the zonal spherical functions which
occur in the principal series of the groups $\PGL_2(\R)$ and
$\PGL_2(\Q_p)$, have been shown to be limits of $q$-ultraspherical
polynomials (see \cite{Askey2} for the $p$-adic limit and
\cite{Koornwinder1} for the real limit). The $p$-adic limit of the
higher rank case appeared in the work of Macdonald \cite{MI1,MI3},
whereas the real limit was proved by Koornwinder\footnote{Lecture at
the INI program on Symmetric Functions and Macdonald Polynomials,
April 2001}. For compact groups, such interpolation has appeared in
the work of Haran \cite{Haran} for the case of the maximal compact
subgroup of $\GL_n$ and its action on the projective line. This has
also been further generalized by Porat \cite{Porat} to invariants of
the $\GL_n$-action on the projective space with respect to upper
triangular matrices.

\subsection{Organization of the paper and notations}
The paper is organized as follows. In section \ref{grassmann} we
describe the Grassmann representation in its various appearances,
the Archimedean in \S\ref{Archimedean}, the non-Archimedean in
\S\ref{non-Archimedean} and the quantum in \S\ref{quantum}. This
section contains a description of all the ingredients required for
carrying out the above plan, with the necessary adjustments and
complements. In section \ref{interpolation}, the ingredients are
glued together to establish the interpolation. Section \ref{example}
contains an example, the one-dimensional case, and section
\ref{finalremarks} is devoted to possible extensions of this work.

{\em Notations:} Whenever possible, the notations of \cite{MI1} have
been followed; partitions are written in a non-increasing order and
are identified with the corresponding Young diagrams. For a
partition $\la=(\la_1,\la_2,\ldots)$, let $\la'$ denote the {\em
transposed} diagram, $|\la|=\sum_i{\la_i}$ its {\em weight} and
$n(\la)=\sum(i-1)\la_i$. The {\em rank} of the partition is the
number of its nonzero parts, and its {\em height} is the largest
part. We shall also use the notation
$\la=(1^{\mu_1}2^{\mu_2}\cdots)$ where $\mu_i=|\{j|\la_j=i\}|$.

Two partial orderings on partitions are used; The partial order
defined by the {\em inclusion of Young diagrams} $\incl$, and the
{\em dominance order} $\domnonstrict$ \footnote{$\la \domnonstrict
\mu \Longleftrightarrow |\la|=|\mu| \mathrm{~and~} \sum_{i=1}^j
\la_i \leq \sum_{i=1}^j \mu_i \quad \forall j \in \N$}. The set of
partitions which consist of at most $m$ parts will be denoted by
$\La_m$. For any ring $A$, we set $\Gt(A)=\GL_n(A)$. In addition to
$q$, three other parameters $(\Al,\Be,\Gam)$ are used. Depending on
the context, they are sometimes rewritten using exponents
$(\al,\be,\gam)$. For $\la=(\la_1,\ldots,\la_m) \in \La_m$, the
vector $(q^{\la_1},\ldots,q^{\la_m})$ is denoted by $q^{\la}$; the
set of all such elements is denoted by $q^{\La_m}$; and
$\rho=(m-1,m-2,\ldots,0)$. $\R$ and $\C$ stand for the real and
complex fields, and $\Kp$ for a non-Archimedean local field with
residue field of cardinality $p^r=|\OO/\wp|$. Multivariable
indeterminants such as $(x_1,\ldots,x_m)$ are abbreviated by $\xb$,
and $\Ac_m=\C[x_1,\ldots,x_m]^{\Sigma_m}$ is the algebra of
symmetric polynomials with $m$ variables. Integration with respect
to any measure, discrete or continuous, is denoted by the integral
sign.

\subsection{Acknowledgements} I am grateful to S. Haran for inspiring discussions,
A. Nevo for hosting this research and commenting on this manuscript,
T. Koornwinder for his hospitality while this work was completed,
 many discussions and careful reading of this manuscript, and to J.
Stokman for suggesting that I look at the shifted Macdonald
polynomials. I thank the referee for his constructive comments.

\section{The Grassmann representation}\label{grassmann}

\subsection{Archimedean theory}\label{Archimedean}
This section is concerned with the Archimedean fields $\R$ and $\C$.
All the objects involved are well known (see James and Constantine
\cite{JC}, and Vretare \cite{Vretare}), but are described here for
completeness. The corresponding maximal compact subgroup $K$ is the
orthogonal group $\Ot(n)$ in the real case, and the unitary group
$\Ut(n)$ in the complex case.

\subsubsection{Space and measure} Points in the space $\Om^{\R}_m=X^{\R}_m \times_{\Ot(n)} X^{\R}_m$
[resp. $\Om^{\C}_m=X^{\C}_m \times_{\Ut(n)} X^{\C}_m$] represent the
relative position of two $m$-dimensional subspaces in the real
[resp. complex] Grassmann manifold modulo the action of $\Ot(n)$
[resp. $\Ut(n)$]. They are given in terms of $m$ critical angles $0
 \le \theta_1 \le \cdots \le \theta_m \le \pi/2$ which are
conveniently rewritten \cite[\S5]{JC} using $u_i=\sin^2(\theta_i)$
to give\footnote{we choose the co-ordinates $\sin^2(\theta_i)$
rather than $\cos^2(\theta_i)$, see \S\ref{nonarchspace}.}
\begin{equation}\label{archsp}
\Om_m^{\R} \simeq \Om_m^{\C} \simeq \Om_m := \{\ub=(u_1,\ldots,u_m)
| 0 \le u_1 \le \cdots \le u_m \le 1\}
\end{equation}
 The projection of the normalized Haar measure from $K$ to the
orbit space $\Om_m$ is given by special values of the parameters in
the Selberg measure \cite{Selberg,AAR} which is given by
\begin{equation}
\dS_m(\ub;\al,\be,\gam)=s_m^{\al,\be,\gam}\prod_{i=1}^{m}u_i^{\al/2-1}(1-u_i)^{\be/2-1}\prod_{i
< j }|u_i-u_j|^{\gam} d\ub
\end{equation} where
\begin{equation}
s_m^{\al,\be,\gam}=\prod_{j=1}^{m}\frac{\Gamma(\al/2+\be/2+(m+j-2)\gam/2)\Gamma(\gam/2)}{\Gamma(\al/2+(j-1)
\gam/2)\Gamma(\be/2+(j-1)\gam/2)\Gamma(j\gam/2)}
\end{equation}
That this is a probability measure on $\Om_m$ is due to Selberg
\cite{Selberg}. We are interested in the following specializations
\begin{align*}
\dh^{\R}(\ub)&=\dS_m(\ub;n-2m+1,1,1)\\
\dh^{\C}(\ub)&=\dS_m(\ub;2(n-2m+1),2,2)
\end{align*}

\subsubsection{Idempotents}

Define an inner product on the algebra of symmetric polynomials
$\Ac_m$ by
\begin{equation}
\langle f,g \rangle_{\al,\be,\gam}=\int_{\Om_m}
f(\ub)\overline{g(\ub)}~\dS_m(\ub;\al,\be,\gam) \qquad f,g \in \Ac_m
\end{equation}
Let $\{\mt_{\la}\}_{\la \in \La_m}$ be the monomial basis of
$\Ac_m$
\begin{equation}
\mt_{\la}(\xb)=\sum_{\eta}x_1^{\eta_1} \cdots x_m^{\eta_m}
\end{equation}
where the summation is over all distinct permutations of $\eta$ of
$\la$. The {\em generalized Jacobi polynomials}
$\{\Eb_{\la}(\xb;\al,\be,\gam)\}_{\la \in \La_m}$ are defined  by
the following conditions \cite{JC,Vretare}
\begin{enumerate}
\item $\Eb_{\la}=d_{\la}\mt_{\la} + \text{lower terms}$, $d_{\la} \neq
0$,
\item $\langle \Eb_{\la},\mt_{\mu} \rangle_{\al,\be,\gam} = 0$  $\forall  \mu
\domstrict \la$,
\item Normalization: $\|\Eb_{\la}\|^2=\Eb_{\la}(\mathbf{0};\al,\be,\gam)$.
\end{enumerate}
Our normalization, which is different from the one in
\cite{JC,Vretare}, is chosen so that the idempotents in the Hecke
algebras are given by the generalized Jacobi polynomials for the
same special values as above:
\begin{align*}
&\eb_{\la}^{\R}(\ub)=\Eb_{\la}(\ub;n-2m+1,1,1) \\
&\eb_{\la}^{\C}(\ub)=\Eb_{\la}(\ub;2(n-2m+1),2,2)
\end{align*}
The generalized Jacobi polynomials are also eigenfunctions with
distinct eigenvalues of a second order differential operator which
specializes to the Laplace-Beltrami operator on the real/complex
Grassmann manifolds after the parameters have been specialized.

\subsection{Non-Archimedean theory}\label{non-Archimedean} Let $\OO$ be
the ring of integers of a non-Archimedean local field $\Kp$. Let
$\wp=(\pi)$ be the maximal ideal and $p^r$ the cardinality of the
residue field $\OO/\wp$. By the principal divisors theorem, any
finite $\OO$-module is of the form $\oplus_{i=1}^{j}\OO/\wp^{\la_i}$
for a partition $\la=(\la_1,\ldots,\la_j)$, which will be referred
to as the {\em type} of the module. As an example, $\La_m$ above
parameterizes all types of finite $\OO$ modules with rank bounded by
$m$. Note that in the partial order defined by the inclusion of
Young diagrams, $\la \le \mu$ if and only if a module of type $\la$
can be embedded in a module of type $\mu$. In such case we shall use
the notation
\begin{equation}
\Big({\mu \atop \la}\Big)= \text{\# of submodules of type $\la$
contained in a module of type $\mu$.}
\end{equation}
Elements in $\La_m$ with height bounded by $k$ are denoted by
$\La_m^{k}=\{\la \in \La_m ~|~ 0 \le \la \le
k^m\}=\{\text{isomorphism types of submodules of
$(\OO/\wp^k)^m$}\}$.

The non-Archimedean theory is completely determined by finite
quotients. More precisely, let $\OO_k=\Ok$ and let $I_k$ stand for
$\Ker\{\Gt(\OO) \rightarrow \Gt(\OO_k)\}$. Each irreducible
representation of the (pro-finite) group $\Gt(\OO)$ factors through
the groups $\Gt(\OO_k)$, except for a finite set of $k \in \N$ whose
cardinality is the {\em level} of the representation. In particular,
the Grassmann representation can be filtered as follows
\begin{equation}\label{filt}\tag{$\star$}
(0) \subset L^2(X^{\Kp}_m)^{I_1} \subset \cdots \subset
L^2(X^{\Kp}_m)^{I_k} \subset \cdots \subset L^2(X^{\Kp}_m)
\end{equation}
and each of its irreducible constituents is contained in some finite
term. The $k$-th term in this filtration is in fact a representation
of $\Gk$, and the direct limit of this sequence is precisely the
smooth part of the Grassmann representation. The finite space $I_k
\backslash X^{\Kp}_m$ can be canonically identified with
$X_{k^m}=\Gr(m,n,\OO_k)$, the Grassmannian of free submodules of
$(\OO_k)^n$ of rank $m$. Thus, we may identify the representation
space $L^2(X^{\Kp}_m)^{I_k}$ with $\FF_{k^m}=\FF(X_{k^m})$, the
space of $\C$-valued functions on $X_{k^m}$.

Let $\HH_{k^m}=\End_{\Gk}(\FF_{k^m})$ be the Hecke algebra
associated with the representation $\FF_{k^m}$. It is isomorphic to
the convolution algebra $\FF(X_{k^m} \times_{\Gk} X_{k^m})$ by
interpreting elements of the latter as $\Gk$-invariant summation
kernels (see \cite[\S2.2]{BO2} for details). The $\Gk$-orbit of an
element $(y,z) \in X_{k^m} \times X_{k^m}$ is determined by the type
of the intersection $y \cap z$, giving rise to the identification
$\La_m^{k} \simeq X_{k^m} \times_{\Gk} X_{k^m}$. The following
diagram summarizes the objects involved and the maps between them.
\[
\begin{matrix}
\Gt(\OO) \quad& X^{\Kp}_m \times X^{\Kp}_m & \rightarrow & ~~X^{\Kp}_m \times_{\Gt(\OO)} X^{\Kp}_m \qquad & ~\dh^{\Kp} \quad& \HH^{\Kp}_m\\
\downarrow \quad & \downarrow& &\downarrow \quad \qquad & \downarrow \quad& \uparrow\\
\Gk  \quad& X_{k^m} \times X_{k^m} & \rightarrow & \La_m^k  \qquad \quad & \dh_k \quad & \HH_{k^m} \\
\downarrow \quad & \downarrow& &\downarrow \quad \qquad & \downarrow \quad& \uparrow\\
\Gt(\OO_{k-1}) \quad& X_{(k-1)^m} \times X_{(k-1)^m} & \rightarrow &
\La_m^{k-1} \qquad& \dh_{k-1} \quad & \HH_{(k-1)^m}
\end{matrix}
\]
Here $\dh^{\Kp}$ and $\dh_k$ stand for the projection of the Haar
measure from $\Gt(\OO)$ to $X^{\Kp}_m \times_{\Gt(\OO)} X^{\Kp}_m$
and $\La_m^k$ respectively. The map from $\La^{k}_m$ to
$\La_m^{k-1}$ is easily described using the transposed Young
diagrams, and is given by $\la'=(\la'_1,\ldots,\la'_{k}) \mapsto
\bar{\la}'=(\la'_1,\ldots,\la'_{k-1})$. We have
\begin{align*}
&\Gt(\OO) \simeq \varprojlim \Gk &&\HH^{\Kp}_m \simeq \varinjlim \HH_{k^m} \\
&X^{\Kp}_m \times_{\Gt(\OO)} X^{\Kp}_m \simeq \hat{\La}_m :=
\varprojlim \La_m^k &&\dh^{\Kp} = \varinjlim \dh_k
\end{align*}

\subsubsection{Space and measure}\label{nonarchspace}

Points in the space $\Om_m^{\Kp} = X_m^{\Kp} \times _{\Gt(\OO)}
X_m^{\Kp}$ correspond to the relative position of two
$m$-dimensional spaces modulo the diagonal $\Gt(\OO)$-action. By the
above discussion it may be identified with $\hat{\La}_m =
\varprojlim \La_m^k$, namely with 'partitions' $(\la_i)_{i=1}^{m}$
where the value $\infty$ is allowed \cite[\S2.3.1]{BO2}. By analogy
with the Archimedean space \eqref{archsp} it is convenient to
rewrite it as
\begin{equation}
\Om_m^{\Kp}=\{\p^{-\la}=(\p^{-\la_1},\cdots,\p^{-\la_m}) | \la \in
\hat{\La}_m \} \subset \Om_m^{\R,\C}
\end{equation}
Note that this embedding is topological, and has the advantage that
the origin $\mathbf{0}=(0,\cdots,0)$ is the common representative of
the orbit of the trivial relative position for all local fields,
namely, $\mathbf{0}=[(x,x)] \in X_m^{\F} \times _{K_{\F}} X_m^{\F}$.
This is also the reason for the choice of the co-ordinates
$u_i=\sin^2(\theta_i)$ rather than $u_i=\cos^2(\theta_i)$ for the
Archimedean spaces.

The following proposition computes the measures $\dh_k$ and
$\dh^{\Kp}$. Note that the measure $\dh^{\Kp}$ vanishes outside the
set $\dot{\Om}_m^{\Kp}=\{\p^{-\la}| \la \in \La_m \} \subset
\Om_m^{\Kp}$.

{\em Notation:} Let $[i]_{q}=1-q^{i}$ for $i \in \N$,
$[i]_q!=[i]_q[i-1]_q \cdots [1]_q$ and $\big[{i \atop
i'}\big]_q=\frac{[i]_q!}{[i']_q![i-i']_q!}$. The index $q$ is
omitted whenever $q=\p^{-r}=|\OO/\wp|^{-1}$.

\begin{prop}\label{padicmeasure}
\begin{align*}
 \mathrm{(1)\quad} &\dh_1(\la)=\frac{\Bigl[{m \atop \la'_1}\Bigr]_{}\Bigl[{n-m \atop m-\la'_1}\Bigr]_{}}{\Bigl[{n
\atop m} \Bigr]_{}}\p^{-r\la'_1(n-2m+\la'_1)} &&  (k=1)\\ & \\
 \textbf{\quad} &\frac{\dh_k(\la)}{\dh_{k-1}(\bar{\la})}
 = \biggl[{\la'_{k-1} \atop \la'_k}\biggr]_{}\frac{[n-2m+\la'_{k-1}]_{}!}{[n-2m+\la'_k]_{}!}
 \p^{-r\la'_k(n-2m+\la'_k)}&&(k>1)\\& \\ {\mathrm{(2)} \quad } &\dh^{\Kp}(\p^{-\la})=\frac{\Bigl[{m \atop
m-\la'_1,\la'_1-\la'_2,\ldots}\Bigr]_{}\frac{[n-m]_{}!}{[m-\la'_1]_{}![n-2m]_{}!}}{\Bigl[{n
\atop m} \Bigr]_{}}\p^{-r\sum {(\la'_i)^2} - r(n-2m)\sum{\la'_i}}
\end{align*}
\end{prop}

\proof The proof of part (2) follows directly from part (1) using
\begin{equation}
\dh^{\Kp}(\p^{-\la})=\prod_{k \ge
1}\frac{\dh_k(\la)}{\dh_{k-1}(\bar{\la})}
\end{equation}
where $\dh_0=1$.

To prove (1), start with $k=1$. The measure $\dh_1$ appeared in
connection with the $q$-Johnson association scheme \cite{DP1}, but
is included here for completeness. Fix spaces $z_1 \subset y_1$ of
dimensions $\la'_1 \le m$ in $(\OO/\wp)^n$. Then
\begin{align*}
&n_1=|\{y | \dim y = m, y \cap y_1 = z_1 \}| = \Big[{n-m \atop
m-\la'_1}\Big]_{\p^r}
\p^{r(m-\la'_1)^2} && \\
&n_2=\text{\# of choices for $z_1$}=\Big({1^{m} \atop 1^{\la'_1}}
\Big)=\Big[{m \atop
\la'_1}\Big]_{\p^r}  && \\
&n_3=\text{\# of $m$-dimensional subspaces in $(\OO/\wp)^n$}=\Big[{n
\atop m}\Big]_{\p^r} &&
\end{align*}
and $\dh_1(\la)={n_1n_2/n_3}$, which together with the relation
$\big[{i \atop i'}\big]_{q}=\big[{i \atop
i'}\big]_{1/q}q^{i'(i-i')}$, gives the desired formula. \\

\noindent For $k>1$, fix two $\OO_k$-modules $f \subset F$ of types
$\f=k^m \le k^n=\Phi$. Let $\bar{F}=F/\wp^{k-1}F$ where $z \mapsto
\bar{z}$ is the quotient map. For any module $y$, let $\ta(y)$
denote its isomorphism type. Then for any $\la \le k^m$
\begin{equation*}\begin{split}
\frac{\dh_k\big(\la\big)}{\dh_{k-1}\big(\bar{\la}\big)}&=\frac{\big({{\Phi}
\atop {\phi}}\big)^{-1}}{\big({\bar{\Phi} \atop
\bar{\phi}}\big)^{-1}}\frac{\big|\big\{z~|~\ta(z \cap f)=\la,
~\ta(z)=\f\big\}\big|}{\big|\big\{\bar{z}~|~\ta(\bar{z} \cap
\bar{f})=\bar{\la}, ~\ta(\bar{z})=\bar{\f}\big\}\big|}
\qquad \text{(Haar $\rightarrow$ counting measure)} \\
&=\frac{\big({\bar{\Phi} \atop \bar{\phi}}\big)}{\big({\Phi \atop
\phi}\big)}\cdot\frac{\big({\phi \atop \la}\big)}{\big({\bar{\phi}
\atop \bar{\la}}\big)}\cdot \frac{\big|\big\{z~|~z \cap f=y_0,
~\ta(z)=\f\big\}\big|}{\big|\big\{\bar{z}~|~\bar{z} \cap
\bar{f}=\bar{y}_0, ~\ta(\bar{z})=\bar{\f}\big\}\big|}
\qquad \text{($y_0$ fixed of type $\la$)} \\
&=\bigg(\p^{-rm(n-m)}\bigg)\bigg(\bigg[{\la'_{k-1} \atop
\la'_{k}}\bigg]_{}\p^{r\la'_k(m-\la'_{k})}\bigg)\bigg(\frac{[n-2m+\la'_{k-1}]_{}!}{[n-2m+\la'_{k}]_{}!}
\p^{r(m-\la'_k)(n-m)}\bigg)
\end{split}\end{equation*}
The computation of the first two terms is straightforward
(alternatively, use the explicit formulas in \cite[\S4.1]{BO2}), and
for the third term we argue as follows. Let $y_0 \subset f$ be a
fixed submodule of type $\la$. Let $\bar{z}$ be a fixed submodule of
$\bar{F}$ of type $\bar{\phi}$ such that
$\bar{f}\cap\bar{z}=\bar{y}_0$. Then $\frac{|\{z|z \cap f=y_0,
~\ta(z)=\f\}|}{|\{\bar{z}|\bar{z} \cap \bar{f}=\bar{y}_0,
~\ta(\bar{z})=\bar{\f}\}|}$ counts the number of submodules $z
\subset F$ of type $\f$ which fit into the following diagram
\[
\begin{matrix}
z & \longmapsto & \bar{z} \\
\cup & & \cup \\
y_0=z \cap f & \longmapsto & \bar{y}_0 \\
\cap &  &   \cap   \\
f & \longmapsto & \bar{f} \\
\end{matrix}
\]
That is, we need to count liftings of $\bar{z}$ which intersect $f$
precisely in $y_0$. First, observe that we may assume $\la'_k=0$.
This amounts to moding out a $k^{\la'_k}$-type summand of $y_0$. A
second observation is that counting different liftings of $\bar{z}$
is equivalent to deforming a fixed basis of a chosen lifting.
Namely, let $z$ be a lifting of $\bar{z}$ and let $\Bc_z$ be a basis
for $z$. Complete this basis to a basis $\Bc_F$ of $F$. If $z'$ is
another lifting of $\bar{z}$, then it has a basis $\Bc_{z'}$ which
is a deformation of $\Bc_z$ with elements from $\wp^{k-1}\Bc_F$.
There are in fact many such bases, however, if we deform only with
elements from $\wp^{k-1}(\Bc_F \setminus \Bc_z)$ we get that $z'=z''
\Longleftrightarrow \Bc_{z'}=\Bc_{z''}$. Putting the last two
observations together, we now fix $z$ which fits into the diagram
above together with a basis $\Bc_z$, and count proper deformations
which also fit the diagram. Let $\Bc_z=\coprod_{i=0}^{k}\Bc_z^i$ and
$\Bc_f=\coprod_{i=0}^{k}\Bc_f^i$ be bases of $z$ and $f$
respectively such that
$\coprod_i\pi^{k-i}\Bc_z^{i}=\coprod_i\pi^{k-i}\Bc_f^{i}$ is a basis
for $y_0$. The assumption $\la'_k=0$ implies that
$\Bc_f^k=\Bc_z^k=\emptyset$. Elements of $\Bc_z \setminus
\Bc_z^{k-1}$ can be deformed arbitrarily, and there are
$|\wp^{k-1}(\Bc_F \setminus \Bc_z)|\cdot|\Bc_z \setminus
\Bc_z^{k-1}|=\p^{r(n-m)(m-\la'_{k-1})}$ such deformations. However,
when deforming the $j$-th basis element of $\Bc_z^{k-1}$, elements
from $\wp^{k-1}(\Bc_f \setminus \Bc_f^{k-1})$ together with the span
of the previously chosen $j-1$ elements must be avoided in order not
to enlarge the intersection $y_0$. Thus, the number of possible
deformations of this element is
$(\p^{r(n-m)}-\p^{r(m-\la'_{k-1}+j-1)})$. Multiplying all
contributions gives the desired result for the third term
\[
(\p^{r(n-m)}-\p^{r(m-\la'_{k-1})})\cdots(\p^{r(n-m)}-\p^{r(m-1)})\cdot
\p^{r(n-m)(m-\la'_{k-1})}
\]
and completes the proof of the proposition. \eproof

\subsubsection{Idempotents}\label{Hecke}

We have the following inner product on $\HH_m^{\Kp}$
\begin{equation}
\langle f,g \rangle_{\Kp}= \int_{\Om_m^{\Kp}} f\bar{g}\dh^{\Kp}
\qquad \forall f,g \in \HH_m^{\Kp}
\end{equation}
The idempotents in the algebra $\HH^{\Kp}_m$, considered as
functions on $\Om_m^{\Kp}$, are orthogonal with respect to the
measure $\dh^{\Kp}$. Since $\dot{\Om}_m^{\Kp}=p^{-\La_m}$ is an open
dense subset of $\Om_m^{\Kp}=p^{-\hat{\La}_m}$ (see
\cite[\S2.3]{BO2}), and carries the full measure of the space by
proposition \ref{padicmeasure}, it suffices to know the restrictions
of functions to $\dot{\Om}_m^{\Kp}$. The explicit computation of the
minimal idempotents in $\HH^{\Kp}_m$ has been carried out in
\cite[\S4.2]{BO2}. The algebra $\HH^{\Kp}_m$ is equipped with the
following natural bases
\begin{description}
\item $\bullet$ $\{\gb^{\Kp}_{\la}\}_{\la \in \La_m}$ - {\em Geometric basis} (delta functions supported on points in
$\dot{\Om}_m^{\Kp}$).
\item $\bullet$ $\{\cb^{\Kp}_{\la}\}_{\la \in \La_m}$ - {\em Cellular basis}.
\item $\bullet$ $\{\eb^{\Kp}_{\la}\}_{\la \in \La_m}$ - {\em Algebraic basis}
(minimal idempotents).
\end{description}
The cellular basis is an intermediate basis which plays an important
role in the non-Archimedean theory and also for the interpolation.
 It is lower triangular with respect to the geometric basis, defined explicitly by
\begin{equation}\label{cellular}
\cb^{\Kp}_{\la}=\sum_{\mu \ge \la}\Big({\mu \atop \la} \Big)
\gb^{\Kp}_{\mu}
\end{equation}
On the other hand it is upper triangular with respect to the
algebraic basis; The subspaces
\[
\JJ^{\Kp}_{\la}=\Span_{\C}\{\cb^{\Kp}_{\mu} | \mu \le \la\} \qquad
\qquad \JJ^{\Kp}_{\la^-}=\Span_{\C}\{\cb^{\Kp}_{\mu} | \mu < \la\}
\]
are in fact ideals, and $\{\JJ^{\Kp}_{\la}/\JJ^{\Kp}_{\la^-}\}_{\la
\in \La_m}$ exhaust the irreducible $\HH_m^{\Kp}$-modules (hence the
term cellular basis). As $\eb^{\Kp}_{\la}$ is by definition the
idempotent which corresponds to the representation
$\JJ^{\Kp}_{\la}/\JJ^{\Kp}_{\la^-}$ we have
\begin{equation}
\langle \eb^{\Kp}_{\la},\cb^{\Kp}_{\mu} \rangle_{\Kp}= 0 \qquad
\forall \mu \le \la
\end{equation}

\subsection{Quantum Grassmannians and some symmetric functions}\label{quantum}
In this section we describe the $q$-objects which interpolate
between the objects related to the local fields. For more details,
the reader is referred to \cite{Kaneko,Aomoto1} for the measure
theoretic considerations (the $q$-Selberg measure), to
\cite{Stokman1,KS} for the spherical functions analysis
(multivariable little $q$-Jacobi polynomials), to
\cite{Okounkov1,Okounkov2,Lassalle1,Lassalle2} for the generalized
binomial coefficients and the shifted Macdonald polynomials and to
\cite{Stokman2} for the description of quantum Grassmannians.
Throughout this section the parameters $q,\Gam,\Al$ and $\Be$ are
used, where the first two are the standard Macdonald parameters. In
some parts restrictions are set on their values.

\subsubsection{The $q$-Selberg measure} \label{qselberg} \quad

The $q$-Selberg measure is a multivariable generalization of the
$q$-beta measure \cite{Askey1}. Let $q,t,\Al,\Be \in (0,1)$ and let
$\rho=(m-1,m-2,\ldots,0)$. Let
\begin{equation}
\Om_m^q=\{q^{\la}t^{\rho}=(q^{\la_1}t^{m-1},q^{\la_2}t^{m-2},\ldots,q^{\la_m})
| \la \in \La_m\} \subset \Om_m^{\R,\C}
\end{equation}
and denote
\begin{align}
&(x)_{\infty}=(x;q)_{\infty} = \prod_{i=0}^{\infty}(1-q^ix) && \\
&(x_1,\ldots,x_l)_{\infty}=\prod_{i=1}^l(x_i)_{\infty} &&
\end{align}
The $q$-Selberg measure is given by
\begin{equation}\label{qselberg}
\begin{split}\dS_m^q(\xb;&\Al,\Be,t)=\\ &\prod_{j=1}^{m}\frac{(\Al t^{m-j},\Be
t^{j-1},t^j,qx_j)_{\infty}}{(\Al \Be t ^{m+j-2},t,q,\Be
x_j)_{\infty}}\Al^{\la_j}t^{2j-2} \prod_{j < i
}\frac{(qx_j/tx_i)_{\infty}}{(tx_j/x_i)_{\infty}}
\big(1-\frac{x_j}{x_i}\big)\end{split}
\end{equation}
for $\xb=q^{\la}t^{\rho} \in \Om_m^q$. Askey conjectured
\cite[\S2]{Askey1} that $\dS_m^q(\xb;\Al,\Be,t)$ is a probability
measure supported on $\Om_m^q$ for $t=q^{\gam}, \gam \in \N$. This
was proved independently by Habsieger \cite{Habsieger} and Kadell
\cite{Kadell1}, and was further generalized by Aomoto
\cite[Proposition 2]{Aomoto1} for any $\gam \in \R_{>0}$. Our
notation follows \cite{Aomoto1} with the dictionary
\begin{align*}
&a \leftrightarrow q^{\alpha-(m-1)(2\gam-1)} && m \leftrightarrow n
\\
&\Be \leftrightarrow
q^{\beta+1}  && x_j \leftrightarrow q^{-1}t_{n-j+1}\\
&t  \leftrightarrow q^{\gamma} && t^{\rho} \leftrightarrow
q^{-1}\xi_F
\end{align*}
and Proposition 2 in \cite{Aomoto1} translates into
\begin{equation}
\int_{\Om_m^q}\dS_m^q= \sum_{q^{\la}t^{\rho} \in
\Om_m^q}\dS_m^q(q^{\la}t^{\rho};\Al,\Be,t)=1
\end{equation}
Note that the order of the variables is reversed with respect to
\cite{Aomoto1,Stokman1,Kaneko} since partitions there are written in
non-decreasing order while here they are written in non-increasing
order. Also, we avoid the use of the $q$-Jackson integral, which is
illuminating when one takes the Archimedean limit but is less
adapted for taking the non-Archimedean limit.

\subsubsection{Multivariable little $q$-Jacobi polynomials}\quad

Define an inner product on the algebra of symmetric polynomials
$\Ac_m$ by
\[
\langle f,g \rangle_{q,\Al,\Be,\Gam}=\int_{\Om_m^q}
f(\omega)\overline{g(\omega)}\dS^q_m(\omega;\Al,\Be,\Gam) \qquad f,g
\in \Ac_m
\]
\begin{defn}[Stokman, \cite{Stokman1}]\label{qjacobi}
The multivariable little $q$-Jacobi polynomials
$\{\Eb_{\la}^q(\xb;\Al,\Be,\Gam)\}_{\la \in \La_m}$ are the unique
polynomials defined by
\begin{enumerate}
\item $\Eb_{\la}^q=d_{\la}\mt_{\la} + \text{lower terms}$, $d_{\la} \neq
0$,
\item $\langle \Eb_{\la},\mt_{\mu} \rangle_{q,\Al,\Be,\Gam} = 0$, $\forall  \mu
\domstrict \la$,
\item Normalization\footnote{note the different normalization comparing to \cite{Stokman1}.}: $\|\Eb_{\la}^{q}
\|^2=\Eb_{\la}^q(\mathbf{0};\Al,\Be,\Gam)$.
\end{enumerate}
\end{defn}
The multivariable little $q$-Jacobi polynomials have interpretation
as zonal spherical functions in the representation of $\Ub_q(n)$
which arises from its action on $\Ub_q(n)/\Ub_q(m) \times
\Ub_q(n-m)$. As we focus only on the zonal spherical functions, the
reader is referred to \cite{Stokman2} for a detailed discussion on
the quantum Grassmannian.

\subsubsection{The basis $\{\Cb_{\la}^q\}$}\quad

In the absence of an explicit formula for the multivariable little
$q$-Jacobi polynomials, a key role in the interpolation between the
idempotents in the Hecke algebras is played by a $q$-version of the
non-Archimedean cellular basis \eqref{cellular}. In short, it
consists of a symmetrized and normalized version of the shifted
Macdonald polynomials. We review their definition and some of their
properties. The only parameters to be used here are $(q,t)$. Let
\begin{align*}
&v_{\la}=v_{\la}(q,t)=\prod_{(i,j) \in
\la}(1-q^{\la_i-j}t^{\la_j'-i+1}) \\
&v'_{\la}=v'_{\la}(q,t)=\prod_{(i,j) \in
\la}(1-q^{\la_i-j+1}t^{\la_j'-i})
\end{align*}
The shifted Macdonald polynomials, also known as interpolation
Macdonald polynomials, were defined in \cite{Okounkov1,Knop,Sahi}.
They were further studied in \cite{Okounkov2}, in which an integral
representation is given and a binomial formula. They are defined as
follows \cite[\S1]{Okounkov2}.

\begin{defn}
The shifted Macdonald polynomials $\{P_{\la}^{\star}(\xb;q,t)\}_{\la
\in \La_m}$ are polynomials in $m$ variables defined by the
following conditions
\begin{enumerate}
\item $P_{\la}^{\star}$ has degree $|\la|$,

\item $P_{\la}^{\star}$ is symmetric in the variables $x_it^{-i}$,

\item\label{ttt} $P_{\la}^{\star}(q^{\mu};q,t)=0$ unless $\la \incl \mu$,

\item $P_{\la}^{\star}(q^{\la};q,t)=(-1)^{|\la|}t^{-2n(\la)}
q^{n(\la')}v'_{\la}$.
\end{enumerate}
\end{defn}

The values of these polynomials on points $q^{\mu}$ with $\la \incl
\mu$ are connected to two-parameter generalized binomial
coefficients, defined in \cite[\S4]{Lassalle1} and
\cite[\S1]{Okounkov2}.

\begin{defn}\label{Lassalle} The generalized binomial
coefficients $\big({\mu \atop \la}\big)_{q,t}$ are defined by the
identity
\[
(v_{\la}')^{-1}P_{\la}(\xb;q,t)\prod_{i=1}^{m}(x_i;q)_{\infty}^{-1}=\sum_{\mu}\Big({\mu
\atop
\la}\Big)_{q,t}t^{n(\mu)-n(\la)}(v_{\mu}')^{-1}P_{\mu}(\xb;q,t)
\]
where the $P_{\la}$'s are the Macdonald polynomials.
\end{defn}
The connection between the generalized binomial coefficients and the
shifted Macdonald polynomials is given by

\begin{equation}\label{specialvalues}
\Big({\mu \atop
\la}\Big)_{q,t}=\frac{P_{\la}^{\star}(q^{\mu})}{P_{\la}^{\star}(q^{\la})},
\qquad \qquad \qquad \qquad \text{\cite[\S7]{Lassalle2}
\cite[\S1]{Okounkov2}}
\end{equation}

The $q$-analogue of the non-Archimedean cellular basis which was
described in \S\ref{Hecke}, is the following symmetrized and
normalized version of the shifted Macdonald polynomials
\begin{defn}\label{q-cellular} The basis $\{\Cb^q_{\la}(\xb;t)\}_{\la \in \La_m}$ of $\Ac_m$
is defined by
\begin{equation}
\Cb^q_{\la}(x_1,\ldots,x_m;t)=\frac{P^{\star}_{\la}(x_1t^{1-m},x_2t^{2-m},\ldots,x_m;q,t)}
{P_{\la}^{\star}(q^{\la};q,t)}~, \qquad \la \in \La_m
\end{equation}
\end{defn}

\section{Interpolation}\label{interpolation}

We are now in a position to state our results concerning the
interpolation. Most of them are multidimensional generalizations of
Haran's work \cite{Haran} regarding interpolation between projective
spaces (Grassmannians of lines) over local fields. By interpolation,
we mean that the $q$-objects described in \S\ref{quantum} have
limits which are the local objects described in \S\ref{Archimedean}
and \S\ref{non-Archimedean}. The functions or measures of which we
take limits are of $m$ variables and might carry one, two or three
parameters, in addition to $q$.

\subsection{Definition of the limits}
Two kinds of limits are considered; For
$f^q=f^q(\xb;\Al,\Be,\Gam)\in \C[\xb]$ define the {\em Archimedean
limit} by
\begin{align}
&[\lim_{\arch}f^q](\ub;\al,\be,\gam)=\lim_{q \rightarrow
1}f^q(\ub;q^{\al/2},q^{\be/2},q^{\gam/2}) && \qquad \qquad \ub \in
\Om_m^{\R,\C}
\end{align}
and, the {\em Non-Archimedean limit} by
\begin{align}
&[\lim_{\nonarch}f^q](\p^{-\la};\al,\be,\gam)=\lim_{q \rightarrow
0}f^q(q^{\la}\p^{-\gam\rho};\p^{-\al}, \p^{-\be},\p^{-\gam}) \qquad
&&  \p^{-\la} \in \Om_m^{\K}
\end{align}
In both limits, the parameter $q$ disappears and the parameters
$\Al,\Be,\Gam$ are replaced by $\al,\be,\gam$. In practice, the
non-Archimedean limit amounts to first substituting
$(\Al,\Be,\Gam)=(\p^{-\al},\p^{-\be},\p^{-\gam})$, and then
substituting $q=0$. To get interpretation of these functions in the
Hecke algebras $\HH_m^{\F}$, set
\begin{equation}\label{parameters}
(\al,\be,\gam)=r(n-2m+1,1,1)
\end{equation}
where for non-Archimedean places $r=[\OO/\wp:\mathbb{F}_p]$, the
degree of the residue field over its prime field, and since this
degree is the same as $[\K:\Q_p]$ for non-ramified extensions, we
set $r=1$ for a real place and $r=2$ for a complex place.

\subsection{Interlude for setting the strategy}
We shall follow the following plan. First, we observe that the
Archimedean and non-Archimedean weak limits of the distribution
\[ f \mapsto \int_{\Om_m^q} f \dS_m^{q} ,\]
are the distributions
\[ f \mapsto \int_{\Om_m^{\F}} f \dh^{\F} ,\]
when substituting the appropriate parameters \eqref{parameters}.
Second, we observe that the flag which is used to define the zonal
spherical functions in the quantum Grassmannian converges to the
flag which is used to define the idempotents $\{\eb_{\la}^{\F}\}$.
We then conclude that the zonal spherical functions in the quantum
Grassmannian converge to the idempotents.

\subsection{Limits of the measure}

\noindent We now prove (the non-Archimedean part of)
\begin{thm}\label{thm-measures}
For any local field $\F$ the measure on the space $X_m^{\F}
\times_{K_{\F}} X_m^{\F}$ is a limit of the $q$-Selberg measure.
\end{thm}

\proof For the Archimedean limit see \cite{Stokman1}; As $q
\rightarrow 1$ the space $\Om_m^q$ approximates the space $\Om_m$,
and the distribution $f \mapsto \int_{\Om_m^q}f\dS_m^q$ weakly
converges to the distribution $f \mapsto \int_{\Om_m}f\dS_m$. In
fact, the possible existence of the Archimedean limit was the main
motivation for introducing the $q$-Selberg measure.

For the non-Archimedean limit, we show that the function $\dh^{\Kp}$
is a limit of the function $\dS_m^{q}$. Substituting a typical
element $\omega_{\la}=q^{\la}t^{\rho}$ in the $q$-Selberg measure
gives
\begin{equation*}
\dS_m^q(q^{\la}t^{\rho};\Al,\Be,\Gam)=f_1 \cdot f_2\cdot f_3~,
\qquad \la \in \La_m
\end{equation*}
where
\begin{align*}
f_1&=\prod_{j=1}^{m}\frac{(\Al t^{m-j},\Be
t^{j-1},t^j)_{\infty}}{(\Al \Be t ^{m+j-2},t,q)_{\infty}}
&&\text{(normalization constant)}\\
f_2&=\prod_{j=1}^m\frac{(q^{\la_j+1}t^{m-j})_{\infty}}
{(\Be q^{\la_j}t^{m-j})_{\infty}}\Al^{\la_j}t^{2\la_j(j-1)}&&\text{(local factors)}\\
f_3&=\prod_{j < i
}\frac{(q^{\la_j-\la_i+1}t^{i-j-1})_{\infty}}{(q^{\la_j-\la_i}t^{i-j+1})_{\infty}}
(1-q^{\la_j-\la_i}t^{i-j})&&\text{(mixed factors)}
\end{align*}
Taking the non-Archimedean limit of these expressions gives
\begin{align*}
&[\lim_{\nonarch}f_1](\p^{-\la};\al,\be,\gam)=
\prod_{j=1}^{m}\frac{(1-\p^{-\al
-(m-j)\gam})(1-\p^{-\be-(j-1)\gam})(1-\p^{-j\gam})}{(1-\p^{-\al-\be -(m+j-2)\gam})(1-\p^{-\gam})}&&\\
&[\lim_{\nonarch}f_2](\p^{-\la};\al,\be,\gam)=\prod_{j=1}^{m}\p^{-\la_j[\al+2(j-1)\gam]}\prod_{\{j:\la_j=0\}}(1-\p^{-\be-(m-j)\gam})^{-1}&&\\
&[\lim_{\nonarch}f_3](\p^{-\la};\al,\be,\gam)= \prod_{\substack{j <
i
\\\la_i=\la_j}} \frac{1-\p^{-(i-j)\gam}}{1-\p^{-(i-j+1)\gam}}
&&
\end{align*}
so that the product of these terms is $[\lim_{\nonarch}
\dS_m^q](p^{-\la};\al,\be,\gam)$. To get the non-Archimedean measure
for the Grassmannian we specialize $(\al,\be,\gam)=r(n-2m+1,1,1)$
and get
\begin{align*}
&[\lim_{\nonarch}f_1](\p^{-\la};r(n-2m+1),r,r)=
\frac{[m]![n-m]!}{[1]^m\bigl[{n \atop
m}\bigr][n-2m]!}&&\\
&[\lim_{\nonarch}f_2](\p^{-\la};r(n-2m+1),r,r)=\frac{1}{[m-\la'_1]!}\p^{-r\sum_{j=1}^{m}{\la_j}(n-2m+2j-1)}&&\\
&[\lim_{\nonarch}f_3](\p^{-\la};r(n-2m+1),r,r)=
\frac{[1]^m}{\prod_{i=0}^k[\la'_i-\la'_{i+1}]!} &&
\end{align*}
which will agree with the second part of proposition
\eqref{padicmeasure} once we show that the exponents of $p$ are the
same, that is
\begin{equation*}
-\sum_{j=1}^{m}{\la_j}(n-2m+2j-1)=-\sum {(\la'_i)^2} -
(n-2m)\sum{\la'_i}
\end{equation*}
However, since $\sum{\la'_j}=\sum{\la_i}$, the last equality reduces
to
\begin{equation*}
\sum_{j=1}^{m}{\la_j}(2j-1)=\sum {(\la'_i)^2}
\end{equation*}
And this equality follows from the fact that both sides evaluate the
cardinality of $\End_{\OO}(\oplus \OO/\wp^{\la_i})$. Thus, we
conclude that
\[
\dh^{\Kp}(\p^{-\la})=[\lim_{\nonarch}
\dS_m^q]\big(p^{-\la};r(n-2m+1),r,r\big). \text{\eproof}\]

\subsection{Limits of functions}
In \cite[\S2]{KO1}, Koornwinder has given an alternative proof for
Haran's non-Archimedean limit of little $q$-Jacobi polynomials which
involves the one variable shifted Macdonald polynomials. This
section consists of a generalization of this proof to the
multidimensional case.

\begin{prop}\label{q-cellular-limit}
\[
[\lim_{\nonarch}\Cb^q_{\la}](\p^{-\mu};\p^{-r})=\cb_{\la}^{\Kp}(\p^{-\mu}).
\]
\end{prop} \proof
Using \eqref{specialvalues}, definition \ref{q-cellular} and the
definition of the non-Archimedean limit, we observe that
\begin{equation}
[\lim_{\nonarch}\Cb^q_{\la}](\p^{-\mu};\p^{-r})=\Cb^q_{\la}(q^{\mu}\p^{-r\rho};\p^{-r})\big|_{q=0}=\Big({\mu
\atop \la}\Big)_{0,\p^{-r}}
\end{equation}
Hence, recalling \eqref{cellular}, we should show that $\big({\mu
\atop \la}\big)_{0,\p^{-r}}=\big({\mu \atop
\la}\big)=\cb^{\Kp}_{\la}(p^{-\mu})$. Indeed,
\begin{equation}\label{Lass}
\Big({\mu \atop
\la}\Big)_{q,t}=t^{-n(\mu)+n(\la)}~J_{\mu/\la}(1,t,t^2,\cdots;q,t)
\qquad \qquad \text{\cite[\S15]{Lassalle1}}
\end{equation}
where $J_{\mu/\la}$ are symmetric functions defined in terms of the
integral Macdonald polynomials $J_{\nu}=v_{\nu}P_{\nu}$ by
\begin{equation}\label{skewJ}
J_{\mu/\la}=\frac{v'_{\mu}}{v'_{\la}}\sum_{\nu}(v'_{\nu})^{-1}f^{\mu}_{\nu,\la}J_{\nu}
\qquad  \text{\cite[\S15]{Lassalle1} and \cite[VI~(7.5)]{MI1}}
\end{equation}
with $f^{\mu}_{\nu,\la}=f^{\mu}_{\nu,\la}(q,t)$ defined by
\begin{equation}
P_{\la}(\xb;q,t)P_{\nu}(\xb;q,t)=\sum_{\mu}f^{\mu}_{\nu,\la}(q,t)P_{\mu}(\xb;q,t)
\qquad \text{\cite[VI~(7.1')]{MI1}}
\end{equation}
We now substitute \eqref{skewJ} into \eqref{Lass} and specialize to
the case $q=0$. As $v'_{\la}(0,t)\equiv1$ we get
\begin{align}
&\Big({\mu \atop
\la}\Big)_{0,t}=t^{-n(\mu)+n(\la)}\sum_{\nu}f^{\mu}_{\nu,\la}(0,t)J_{\nu}(1,t,t^2,\cdots;0,t)
\end{align}
however,
\begin{align*}
&J_{\nu}(1,t,t^2,\cdots;q,t)=t^{n(\nu)}  && \text{\cite[p.
366(9)]{MI1}}\\
&f^{\mu}_{\nu,\la}(0,t)=t^{n(\mu)-n(\la)-n(\nu)}g^{\mu}_{\nu,\la}(t^{-1})
&& \text{\cite[p. 217(3.6) and p. 343(7.2)(ii)]{MI1}}
\end{align*}
where $g^{\mu}_{\nu,\la}$ are the Hall polynomials \cite[chapter
II]{MI1}. Thus
\begin{align}
\Big({\mu \atop \la}\Big)_{0,t}&=\sum_{\nu}g^{\mu}_{\nu,\la}(t^{-1})
\end{align}
Since $g^{\mu}_{\nu,\la}(\p^r)$ is by definition the number of
$\OO$-submodules of type $\la$ and co-type $\nu$ in an $\OO$-module
of type $\mu$, summing over all the co-types $\nu$, gives the total
number of submodules of type $\la$, and we have
\[
[\lim_{\nonarch}\Cb^q_{\la}](\p^{-\mu};\p^{-r})=\Big({\mu \atop
\la}\Big)_{0,\p^{-r}}=\Big({\mu \atop
\la}\Big)=\cb^{\Kp}_{\la}(p^{-\mu}). \text{\eproof}
\]

\noindent With this in hand we can prove (the non-Archimedean part
of)
\begin{thm}\label{thm-idempotents}
For any local field $\F$ the idempotents in the Hecke algebra
associated with the Grassmann representation are limits of
multivariable little $q$-Jacobi polynomials.
\end{thm}

\proof The partial orderings $\incl$ and $\domnonstrict$ can be
completed simultaneously to a total ordering,
 e.g. the lexicographical ordering. Let $\mathcal{M}=\text{Flag}\{\mt_{\la}|\la \in \La_m \}$ be the flag
defined by the monomial basis of $\Ac_m$ with respect to such total
ordering. The multivariable little $q$-Jacobi are obtained by
applying the Gram-Schmidt procedure to the flag $\mathcal{M}$, with
respect to the inner product $\langle \cdot, \cdot
\rangle_{q,\Al,\Be,\Gam}$ (definition \ref{qjacobi}).

The Archimedean limit, $\Eb_{\la}=\lim_{\arch}\Eb_{\la}^q$, follows
as this inner product deforms continuously to the inner product
$\langle \cdot, \cdot \rangle_{\al,\be,\gam}$, which is used to
define the generalized multivariable polynomials
(\S\ref{Archimedean}), see \cite{KS} for details.

As for the non-Archimedean limit, we observe that the flag
$\mathcal{M}$ is also defined by the basis $\{\Cb_{\la}^q\}_{\la \in
\La_m}$; Indeed,
\begin{align*}
&\Cb_{\la}^q = P^{\star}_{\la}(q^{\la})^{-1}P_{\la}+\text{lower
terms w.r.t. $\incl$} &&
\text{(by the binomial formula \cite[(1.12)]{Okounkov2}}) \\
&P_{\la}=\mt_{\la}+\text{lower terms w.r.t. $\domnonstrict$} &&
\text{(by definition \cite[VI(4.7)]{MI1})}
\end{align*}
Thus, using the total ordering which refines both partial orderings,
$\{\Cb_{\la}^q\}$ and $\{\mt_{\la}\}$ define the same flag
$\mathcal{M}$. The idempotents basis in the non-Archimedean Hecke
algebra $\HH^{\Kp}_m$ are obtained by applying the Gram-Schmidt
procedure to the cellular basis $\{\cb^{\Kp}_{\la}\}$. As the
$q$-Selberg measure deforms continuously to the non-Archimedean
measure $\dh^{\Kp}$ (by theorem \ref{thm-measures}), and the basis
$\{\Cb_{\la}^q\}$ converges to the cellular basis (by proposition
\ref{q-cellular-limit}), the multivariable little $q$-Jacobi
polynomials converge to the idempotents up to constants. Our
normalization in definition \ref{qjacobi} is designed to eliminate
these constants, as for idempotents one has
$\|\eb^{\Kp}_{\la}\|^2=\eb^{\Kp}_{\la}(\mathbf{0})$. \eproof

\section{Example}\label{example}

This section is devoted to the one-dimensional case which was
treated in \cite{Haran} (see also \cite{KO1}), as it admits a
completely explicit description. For m=1, the representation of
$K^{\F}$ arises from its action on the projective space
$X_1^{\F}=\PP^{n-1}_{\F}$. The representation $L^2(\PP^{n-1}_{\F})$
decomposes into irreducible representations
$\{\mathcal{U}^{\F}_{\la}\}_{\la \in \La_1}$ where $\La_1=\N_0$. The
space $\PP^{n-1}_{\F} \times_{K^{\F}}\PP^{n-1}_{\F}$, which
describes the $K^{\F}$-relative position of two lines, is given by
$[0,1]$ (normalized angles) for an Archimedean place and by
$\{\p^{-\la}\}_{\la \in \N_0\cup\{\infty\}} \subseteq [0,1]$ for a
non-Archimedean place. The triplets (space, measure, idempotents)
are given as follows.

\begin{description}

\item \textbf{Archimedean.} For $u \in [0,1]$ and $\la \in \N_0$ let
\begin{align*}
&\dS(u;\al,\be)=\genfrac{}{}{0.1pt}{1}{\Gamma(\frac{\al}{2}+\frac{\be}{2})}{\Gamma(\frac{\al}{2})\Gamma(
\frac{\be}{2})}
u^{\frac{\al}{2}-1}(1-u)^{\frac{\be}{2}-1}du&  && \\
&\Eb_{\la}(u;\al,\be)=
\genfrac{}{}{0.1pt}{1}{(\frac{\al}{2})_{\la}(\frac{\al}{2}+\frac{\be}{2})_{\la}}
{(\frac{\be}{2})_{\la}~\la!}\genfrac{}{}{0.1pt}{1}{2\la-1+\frac{\al}{2}+\frac{\be}{2}}{\la-1+\frac{\al}{2}+\frac{\be}{2}}
  \mbox{}_2F_1\big[
\genfrac{}{}{0pt}{1}{-\la,\la+\frac{\al}{2}+\frac{\be}{2}-1}{\al};u\big]&
&&
\end{align*}
where $(y)_j=y(y+1)\cdots(y+j-1)$ is the shifted factorial and
$\mbox{}_2F_1$ is the hypergeometric
function\footnote{$\mbox{}_2F_1\Big[
\genfrac{}{}{0pt}{0}{\alpha_1,\alpha_2}{\alpha_3};u\Big]=\sum_{j=0}^{\infty}\frac{(\alpha_1)_j(\alpha_2)_j}{(\alpha_3)_jj!}u^j
$}. $\dS(u;\al,\be)$ is the normalized beta measure on the unit
interval, and $\{\Eb_{\la}(u;\al,\be)\}_{\la \in \N_0}$ are the
normalized Jacobi polynomials. For the special values
$(\al,\be)=(n-1,1)$ and $(\al,\be)=2(n-1,1)$ the triplet
$([0,1],\dS(u;\al,\be),\{\Eb_{\la}(u;\al,\be)\}_{\la \in \N_0})$
specializes to the real and complex triplets
$([0,1],\dh^{\R},\{\eb^{\R}_{\la}\}_{\la \in \N_0})$ and
$([0,1],\dh^{\C},\{\eb^{\C}_{\la}\}_{\la \in \N_0})$. The dimensions
of the irreducible representations are given for by
\begin{align*}
&\dim
\Uc_{\la}^{\R}=\Eb_{\la}(0;n-1,1)=\frac{2\la+\frac{n}{2}-1}{\la+\frac{n}{2}-1}\frac{(\frac{n-1}{2})_{\la}(\frac{n}{2})_{\la}}
{(\frac{1}{2})_{\la}\la!}~, && \qquad \\
&\dim
\Uc_{\la}^{\C}=\Eb_{\la}(0;2n-2,2)=\frac{2\la+n-1}{n-1}\Big({n+\la-2
\atop \la}\Big)^2~. && \qquad
\end{align*}

\item \textbf{Non-Archimedean.} The case $m=1$ is greatly simplified by the fact that the terms in the
filtration \S\ref{non-Archimedean}\eqref{filt} are in bijection with
the irreducibles, and each step in the filtration contains exactly
one new irreducible representation. It follows that $\la$ and $k$
are identified and
\begin{align*}
\dim
\Uc_{\la}^{\K}=|\PP^{n-1}_{\OO/\wp^{\la}}|-|\PP^{n-1}_{\OO/\wp^{\la-1}}|=\left\{
\begin{array}{ll}
    \frac{(1-p^{-r(n-1)})}{(1-p^{-r})}p^{r(n-1)}, & \la=1 \\
    \frac{(1-p^{-rn})(1-p^{-r(n-1)})}{(1-p^{-r})}p^{r(n-1)\la}, & \la \ge 2 \\
\end{array}
\right.
\end{align*}
and is equal to $1$ for $\la=0$, where $|\PP^{n-1}_{\OO/\wp^{\la}}|=
\frac{1-\p^{-rn}}{1-\p^{-r}}p^{r(n-1)\la}$ for $\la \ge 1$, and
$|\PP^{n-1}_{\OO/\wp^{0}}|=1$. The measure is easily seen to be
\begin{align*}
&\dh^{\Kp}(\p^{-\la})=\left\{
\begin{array}{ll}
    \genfrac{}{}{0.1pt}{1}{|\PP^{n-1}_{\OO/\wp}|-1}{|\PP^{n-1}_{\OO/\wp}|}=\frac{1-\p^{-r(n-1)}}{1-\p^{-rn}}, &\quad \la=0 \\
    \genfrac{}{}{0.1pt}{1}{1}{|\PP^{n-1}_{\OO/\wp}|}\genfrac{}{}{0.1pt}{1}{|\Abb^{n}_{\wp/\wp^{\la}}|-1}
    {|\Abb^{n}_{\wp/\wp^{\la}}|}=\genfrac{}{}{0.1pt}{1}{(1-\p^{-r})(1-\p^{-r(n-1)})}{(1-\p^{-rn})\p^{-r(n-1)\la}}, &\quad \la \ge 1  \end{array}
\right.
\end{align*}
where $\Abb^n$ stands for the affine $n$-space. For the idempotents,
we use again the fact that the filtration admits only one new
irreducible in each step, but this time on the level of the Hecke
algebras. The Hecke algebra $\HH_1^{\K}$ is the direct limit the
algebras $\{{\HH}_{\la}\}_{\la \in \N_0}$ ($\la=k^1$ of
\S\ref{non-Archimedean}). Each of these algebras contains a unit
element $\mathbf{1}_{\la}$, which as a function on the orbits space
is given by $|\PP^{n-1}_{\OO/\wp^{\la}}|{\bbbone}_{\{\p^{-\mu}| \mu
\ge \la\}}$. Thus on $p^{-\N_0}$ we have
\begin{align*}
&\eb^{\Kp}_{\la}=\left\{
\begin{array}{ll}
    \mathbf{1}_0={\bbbone}_{\p^{-\N_0}}, & \la=0 \\
    \mathbf{1}_{\la}-\mathbf{1}_{\la-1}=|\PP^{n-1}_{\OO/\wp^{\la}}|{\bbbone}_{\{\p^{-\mu}| \mu
\ge \la\}}-|\PP^{n-1}_{\OO/\wp^{\la-1}}|{\bbbone}_{\{\p^{-\mu}| \mu
\ge \la-1\}}, &\la \ge 1 \\
\end{array} \right. \\ \\
&\cb^{\K}_{\la}=\sum_{\mu \ge
\la}\gb^{\K}_{\mu}={\bbbone}_{\{\p^{-\mu}|\mu \ge \la\}}
\end{align*}
and the non-Archimedean triplet is $(\p^{-\N_0} \cup
\{0\},\dh^{\Kp}, \{ \eb^{\Kp}_{\la} \}_{\la \in \N_0})$.

\item \textbf{Quantum.} For $q \in (0,1)$ and $\la \in \N$ let
\begin{align*}
&\dS^q(q^{\la};\Al,\Be)=\frac{(\Al;q)_{\infty}}
{(\Al\Be;q)_{\infty}}\frac{(\Be;q)_{\la}}
{(q;q)_{\la}}\Al^{\la} &&  \\
&\Eb^q_{\la}(x;\Al,\Be)= \frac{(1-\Al\Be q^{2\la-1})(\Al\Be
q^{-1};q)_{\la}(\Al;q)_{\la}}{(1-\Al\Be
q^{-1})(q;q)_{\la}(\Be;q)_{\la}\Al^{\la}}\mbox{}_2\phi_1\Big[{q^{-\la},q^{\la-1}
\Al \Be \atop \Al };q,qx\Big] && \\
&\Cb^q_{\la}(x)=\frac{(x;q^{-1})_{\la}}{(q^{\la};q^{-1})_{\la}}=\frac{(x-1)(x-q)\cdots(x-q^{\la-1})}{(q^{\la}-1)(q^{\la}-q)\cdots(q^{\la}-q^{\la-1})}&&
\end{align*}
where $(y;q)_j=(1-y)(1-y q)\cdots(1-y q^{j-1})$ is the $q$-shifted
factorial and $\mbox{}_2\phi_1$ is the basic hypergeometric
function\footnote{$\mbox{}_2\phi_1\Big[
\genfrac{}{}{0pt}{0}{a_1,a_2}{a_3};u\Big]=\sum_{j=0}^{\infty}\frac{(a_1;q)_j(a_2;q)_j}{(a_3;q)_j(q;q)!}u^j$.
}. $\dS^q(q^{\la};\Al,\Be)$ is the normalized $q$-beta measure on
the set $\{q^{\la}\}_{\la=0}^{\infty}$ and the associated orthogonal
base consists of the normalized little $q$-Jacobi polynomials
\cite{GR,KO1}, $\{\Eb^q_{\la}(x;\Al,\Be)\}_{\la \in \N_0}$. Then the
$q$-triplet is given by $(q^{\N_0}, ~\dS^q(q^{\la};\Al,\Be),
\{\Eb^q_{\la}(x;\Al,\Be)\}_{\la \in \N_0})$.
\end{description}

{\em Remarks}
\begin{itemize}
\item The parameter $t$ does not appear in the one dimensional case.
\item The formula
\[
\texttt{D}_{\la}^q(\Al,\Be)=\Eb^q_{\la}(\mathbf{0};\Al,\Be)=\frac{(1-\Al\Be
q^{2\la-1})(\Al\Be q^{-1};q)_{\la}(\Al;q)_{\la}}{(1-\Al\Be
q^{-1})(q;q)_{\la}(\Be;q)_{\la}\Al^{\la}}
\]
interpolates between the dimensions of the irreducible
representations $\Uc_{\la}^{\F}$.
\item The non-Archimedean limit for $\dS^q$, $\texttt{D}_{\la}^q$ and $\Cb^q_{\la}$ is
immediate. For more details regarding this limit see \cite{Haran}
and \cite[\S2]{KO1}. The Archimedean limit of the $q$-beta measure
and basic hypergeometric series is discussed in \cite[\S1]{Askey1}
and \cite[pp. 1-28]{GR}.
\end{itemize}

\section{Related problems} \label{finalremarks}

\subsection{The module of intertwining operators
$\Sc(X_{m_1}^{\F} \times_{K^{\F}} X_{m_2}^{\F})$} For $m_1 \le m_2
\le [n/2]$, one can consider in a similar manner the module of
intertwining operators between the representations
$\Sc(X_{m_2}^{\F})$ and $\Sc(X_{m_1}^{\F})$. This results in a very
similar discussion, where the only difference occurs in the
parameters $\al,\be$ and $\gam$, while the geometry remains as in
the equal dimension case for $m=m_1$. As an example see \cite{KO1}
for the case $m_1=1$.

\subsection{Dimensions of the irreducible representations}
The $q$-dimension of the irreducible representation $\Uc_{\la}^{q}$
in the quantum Grassmannian, which independently on the
normalization is given by
\begin{equation}
\texttt{D}_{\la}^q(\Al,\Be,\Gam)=\frac{\Eb^q_{\la}(\mathbf{0};\Al,\Be,\Gam)^2}{\|\Eb^q_{\la}(\xb;\Al,\Be,\Gam)\|^2}
\end{equation}
interpolates between the dimensions of the irreducible
representations which correspond to $\la$ for all local fields.

\subsection{Haran's process} The case $m=1$ was studied extensively by Haran in
\cite{Haran}. Haran also constructs discrete random processes in
order to obtain the bases for the Archimedean places, the
non-Archimedean places and the $q$-case. The bases are defined on
the Martin boundary of the processes. It would be interesting to
find a generalization of these processes in the case of
Grassmannians.

\subsection{Other algebraic groups} A natural venue for further study is to consider other
maximal compact subgroups $K^{\F}$ of reductive algebraic groups and
natural multiplicity free representations of them $V^{\F}$. The
finite analogue of such representations can be found in \cite{SD1},
in which such representations of Chevalley groups over finite fields
are studied. These can be considered as the \emph{level zero} part
of representations of the maximal compact subgroups. Roughly, the
picture is
\begin{align*}
V^{\R}, V^{\C} \longleftarrow q \longrightarrow & ~V^{\Kp} \\
             &  ~\uparrow  \\
             &             ~V^{\OO/\wp}=\text{level zero part of}~V^{\Kp}
\end{align*}
For example, the particular case of $\GL_n(\OO/\wp)$, which admits
the $|\OO/\wp|^{-1}$-Hahn polynomials as idempotents, is just the
first term in the filtration described in \S\ref{non-Archimedean}.







\end{document}